\documentclass[12pt,a4paper,notitlepage]{article}
\usepackage{amsmath,amsfonts,amsthm,amssymb}
\usepackage{mathpple}
\newcommand\bm{{\mathbf m}} 
\newcommand\bn{{\mathbf n}}

\newcommand\bv{{\mathbf v}}

\newcommand\bzero{{\mathbf 0}}

\newcommand\bA{{\mathbf A}}

\newcommand\bJ{{\mathbf J}} 
 
\newcommand\bL{{\mathbf L}} 
\newcommand\bM{{\mathbf M}} 
\newcommand\bN{{\mathbf N}}

\newcommand\bR{{\mathbf R}}

\newcommand\bV{{\mathbf V}} 
\newcommand\bW{{\mathbf W}} 
\newcommand\bX{{\mathbf X}}

\newcommand\efg{{\mathfrak g}}

\newcommand\cD{{\mathcal D}}

\newcommand\cG{{\mathcal G}}

\newcommand\cJ{{\mathcal J}}

\newcommand\cM{{\mathcal M}} 
\newcommand\cN{{\mathcal N}}
\newcommand\cO{{\mathcal O}}  
\newcommand\cP{{\mathcal P}}

\newcommand\cT{{\mathcal T}}

\newcommand\mvector{\boldsymbol}

\newcommand\vv{\mvector{v}}
\newcommand\vw{\mvector{w}}

\newcommand\field{\mathbb}
\newcommand\R{\field{R}}
\newcommand\bbS{\mathbb{S}}
\newcommand\C{\field{C}}
\newcommand\Z{\field{Z}}
\newcommand\N{\field{N}}
\newcommand\Q{\field{Q}}
\newcommand\bbP{\mathbb{P}}

\newcommand\bOmega{\boldsymbol{\Omega}}

\renewcommand\Re{\operatorname{Re}}

\newcommand\ord{\operatorname{ord}}

\newcommand\rmd{\mathrm{d}}
\newcommand\SLtwoC{$\mathrm{SL}(2,\C)$}
\newcommand\CPOne{$\C\mathbb{P}^1$}
\newcommand\rmi{\mathrm{i}}

\newcommand\Dt{\frac{\rmd\phantom{t} }{\rmd t}} 
\newcommand\Dtt{\frac{\rmd^2\phantom{t} }{\rmd t^2}} 
 
\newcommand\rscalar[2]{\langle #1, #2\rangle} 
\newcommand\cn{\operatorname{cn}}
\newcommand\sn{\operatorname{sn}}
\newcommand\dn{\operatorname{dn}}

\usepackage[square]{natbib}
\usepackage{graphicx}
\theoremstyle{plain}
\newtheorem{theorem}{Theorem}
\newtheorem{lemma}{Lemma}
\newtheorem{proposition}{Proposition}
\theoremstyle{definition}

\newtheorem{remark}{Remark}
\begin{document}
\thispagestyle{plain}
\vspace{2em}
\hfill\textsf{accepted to: Annales de la Faculte des Sciences de   Toulouse}
\vspace{2em}
\begin{center}
{\Large
\textbf{Differential Galois Approach to the Non-integrability of the
  Heavy Top Problem}}\\[3em]
{\large\textbf{Andrzej J.~Maciejewski}}\\[1em]
  Institute of Astronomy,
  University of Zielona G\'ora, \\
  Podg\'orna 50, PL-65--246 Zielona G\'ora, Poland,
  (e-mail: maciejka@astro.ia.uz.zgora.pl)\\[1em]
{\large\textbf{  Maria Przybylska} }\\[1em] 
  INRIA Projet \textsc{Caf\'e},\\
  2004, Route des Lucioles, B.~P. 93,
  06902 Sophia Antipolis Cedex, France,  \\
  and \\
Toru\'n Centre for Astronomy,
  Nicholaus Copernicus University, \\
  Gagarina 11, PL-87--100 Toru\'n, Poland,
  (e-mail: Maria.Przybylska@sophia.inria.fr) 
\end{center} 

\vspace{2em}

\begin{flushright}
   \textsf{---To Jean-Pierre Ramis} 
  \end{flushright}
  \vspace{2em}
  {\small \textbf{Abstract}\\
    We study integrability of the Euler-Poisson equations describing
    the motion of a rigid body with one fixed point in a constant
    gravity field. Using the Morales-Ramis theory and tools of
    differential algebra we prove that a symmetric heavy top is
    integrable only in the classical cases of Euler, Lagrange, and
    Kovalevskaya and is partially integrable only in the
    Goryachev-Chaplygin case. Our proof is alternative to that given
    by Ziglin ({\em Funktsional. Anal. i Prilozhen.}, 17(1):8--23,
    1983; {\em Funktsional. Anal. i Prilozhen.}, 31(1):3--11, 95, 1997).}
  \\[1em]
  {\small \textbf{R\'esum\'e}\\
    Nous \'etudions l'int\'egrabilit\'e des \'equations de
    Euler-Poisson qui d\'ecrivent le mouvement d'un solide rigide avec
    un point fixe dans un champ gravitationnel constant. En utilisant
    la th\'eorie de Morales-Ramis et des outils d'alg\`ebre
    diff\'erentielle, nous prouvons qu'un solide sym\'etrique est
    int\'egrable seulement dans les cas classiques d'Euler, Lagrange
    et Kowalevski, et est partiellement int\'egrable seulement dans le
    cas Goryatchev-Tchaplygin. Notre preuve est une alternative \`a
    celle donn\'ee par Ziglin ({\em Funktsional. Anal. i Prilozhen.},
    17(1):8--23, 1983; {\em Funktsional. Anal. i Prilozhen.},
    31(1):3--11, 95, 1997).}  \newpage

%
%

\section{Equations of motion and motivation}
\label{sec:eqandm}
Equations of motion of a rigid body in external fields are usually
written in a body fixed frame. Here, we use the following convention.
For a vector $\vv$ we denote by $\bV=[V_1,V_2,V_3]^T$ its coordinates
in a body fixed frame, and we consider it as a one column matrix.  The
vector and scalar products of two vectors $\vv$ and $\vw$ expressed in
terms of the body fixed coordinates are denoted by $[\bV,\bW]$ and
$\rscalar{\bV}{\bW}$, respectively.

We consider a rigid body of mass $m$ located in a constant gravity
field of intensity $g$. One point of the body is fixed. The distance
between the fixed point and the mass centre of the body is $D$. Assuming
that $gD\neq 0$, we choose units in such a way that $\mu:=mgD=1$. The
Euler-Poisson equations
\begin{equation}
  \label{eq:ep}
  \Dt \bM =[ \bM,\bJ\bM]+ [\bN, \bL], \qquad
  \Dt \bN = [\bN , \bJ\bM],  
\end{equation}
describe the rotational motion of the body. In the above equations
$\bM$ denotes the angular momentum of the body, $\bN$ is the unit
in the direction of the gravity field, and $\bL$ is the
unit vector from the centre of mass of the body to the fixed point;
$\bJ$ is the inverse of the matrix of inertia, so $\bOmega:=\bJ\bM$ is
the angular velocity of the body. The principal moments of inertia
$A$, $B$, and $C$ are eigenvalues of $\bJ^{-1}$.  For our further
consideration it is important to notice that in~\eqref{eq:ep} the body
fixed frame is unspecified, so we can choose it according to our
needs.  A body fixed frame in which $\bJ$ is diagonal is called the
principal axes frame. This frame is uniquely defined (up to the
numbering of the axes) when $\bJ$ has no multiple eigenvalues. If
$\bJ$ has a multiple eigenvalue, e.g. when $A=B$, we say that the body
is symmetric. Then the principal axes frame is defined up to a
rotation around the symmetry axis.

System~\eqref{eq:ep} depends on parameters $A$, $B$, $C$, and $\bL$
but physical constraints restrict the allowable values of parameters to
a set $\cP\subset\R^6$ defined by  the following conditions 
\begin{gather*}
\label{eq:par}
A>0,\quad B>0, \quad C>0, \qquad \rscalar{\bL}{\bL}=1, \\
A+B\geq C, \quad B+C\geq A, \quad C + A\geq B.
\end{gather*}

Euler-Poisson equations possess three first integrals
\begin{gather}
\label{eq:H}
H = \frac{1}{2}\rscalar{\bM}{\bJ\bM} +  \rscalar{\bN}{\bL},\\
\label{eq:H1H2}
H_1 = \rscalar{\bM}{\bN}, \qquad H_2 = \rscalar{\bN}{\bN}.
\end{gather}
It is known that  on the level
\begin{equation*}
  \label{eq:m4}
   M_\chi = \{ (\bM,\bN)\in \R^6\  | \ H_1 = \chi, H_2 = 1 \}
\end{equation*}
the Euler-Poisson equations are the  Hamiltonian ones, see e.g.
\cite{Arnold:78::,Marsden:94::,Audin:96::}.  
\begin{remark}
\label{rem:r1}
\textsl{ Configuration space of a rigid body with a fixed point is the
  Lie group $\mathrm{SO}(3,\R)$ (all possible orientations of the body
  with respect to an inertial frame).  Thus, classically, the phase
  space for the problem is $T^*\mathrm{SO}(3,\R)$.  Hence it is a
  Hamiltonian system with three degrees of freedom and it possesses
  one additional first integral $H_1$ (first integral $H_2$ is
  identically equal to one in this formulation). The existence of this
  first integral is related to the symmetry of the problem (rotations
  around the direction of the gravity field) and this allows to reduce
  the system by one degree of freedom.  The Euler-Poisson equations
  can be viewed as an effect of reduction of the system on
  $T^*\mathrm{SO}(3,\R)$ with respect to this symmetry.  The phase
  space of the reduced problem can be considered as the dual $\efg^*$
  to Lie algebra $\efg$ of group of rigid motions $G=\R^3
  \rtimes\mathrm{SO}(3,\R)$. We can identify $(\bN,\bM)$ as element of
  $\efg^*={\R^{3*}}\rtimes \mathrm{so}(3,\R)^*$ using standard
  isomorphisms between $ \mathrm{so}(3,\R)$ and $\R^{3}$, and between
  ${\R^{3*}}$ and $\R^3$. Let $(\bX,\bA)\in G$, with $\bX\in\R^3$ and
  $\bA \in \mathrm{SO}(3,\R)$.  Then the coadjoint action
  $\mathrm{Ad}^*_{(\bX,\bA)}:\efg^*\rightarrow\efg^* $ is defined by
\[
\mathrm{Ad}^*_{(\bX,\bA)}(\bN,\bM)=( \bA\bN,[\bX,\bA \bN] + \bA\bM).
\]
Functions $H_1$ and $H_2$ given by~\eqref{eq:H1H2} are orbits
invariant, i.e., on each orbit of coadjoint action 
\[
\cO_{(\bN,\bM)}:= \bigcup_{(\bX,\bA)\in G}\{( \bA\bN,[\bX,\bA \bN] + \bA\bM)\}
\]
they have constant
values. They are called Casimirs. Orbits of coadjoint action defined
above coincide with $M_\chi$ which is diffeomorphic to $T\bbS^2$. As
it is well known, orbits of a coadjoint action are symplectic
manifolds equipped with the standard Kostant-Berezin-Souriau-Kirillov
symplectic structure \cite{Kirillov:76::,Souriau:70::}.  On four
dimensional orbits $M_\chi$ the Euler-Poisson equations are
Hamiltonian with $H$ given by~\eqref{eq:H} as the Hamilton
function.  }
\end{remark}
Thus, as a Hamiltonian system on $M_\chi$ the Euler-Poisson equations
have two degrees of freedom and are integrable on $M_\chi$ if there
exists an additional first integral $H_3$ which is functionally
independent with $H$ on $M_\chi$.  Equivalently, we say that the
Euler-Poisson equations are integrable if there exists a first
integral $H_3$ which is functionally independent together with $H$,
$H_1$ and $H_2$.

We say that the Euler-Poisson equations are partially integrable if
they are integrable on $M_0$.

The known integrable cases are the following:
\begin{enumerate}
\item The Euler case (1758) corresponds to the situation when there is
  no gravity (i.e. when $\mu=0$) or $\bL=\bzero$ (the fixed point of
  the body is the centre of mass). The additional first integral in
  this case is the total angular momentum $H_3 = \rscalar{\bM}{\bM}$.
\item In the Lagrange case~\cite{Lagrange:1889::} the body is
  symmetric (i.e. two of its principal moments of inertial are equal)
  and the fixed point lies on the symmetry axis. The additional first
  integral in this case is the projection of the angular momentum onto
  the symmetry axis. If we assume that $A=B$, then in the Lagrange
  case $L_1=L_2=0$, and $H_3=M_3$.
\item In the Kovalevskaya
  case~\cite{Kowalevski:1888::,Kowalevski:1890::} the body is
  symmetric and the principal moment of inertia along the symmetry
  axis is half of the principal moment of inertia with respect to an
  axis perpendicular to the symmetry axis. Moreover, the fixed point
  lies in the principal plane perpendicular to the symmetry axis. If
  $A=B=2C$, then (after an appropriate rotation around the symmetry
  axis) we have in the Kovalevskaya case $L_2=L_3=0$. The additional
  first integral has the form
 \[
H_3 = \left( \frac{1}{2}(M_1^2-M_2^2)+N_1\right)^2 +(M_1M_2+N_2)^2.
\] 
\item In the Goryachev-Chaplygin case~\cite{ Goryachev:1910::} the body
  is symmetric and, as in the Kovalevskaya case, the fixed point lies
  in the principal plane perpendicular to the symmetry axis. If we
  assume that the third principal axis is the symmetry axis, then in
  the Goryachev-Chaplygin case we have $A=B=4C$ and $L_2=L_3=0$. In
  the Goryachev-Chaplygin case equations~\eqref{eq:ep} are integrable
  only on the level $H_1=0$ and the additional first integral has the
  following form:
 \[
H_3 =
  M_3(M_1^2+M_2^2) + M_1N_3.
\]
\end{enumerate} 
For a long time the question if there are other integrable cases of
the Euler-Poisson equations except those enumerated above was open,
although many leading mathematicians tried to give a positive answer
to it.

The problem was completely solved by S.L.~Ziglin in a brilliant way.
First, in \cite{Ziglin:80::d} he proved the following.
\begin{theorem}[Ziglin, 1980]
  If $(A-B)(B- C)(C-A)\neq 0$, then the Euler-Poisson system
  does not admit a real meromorphic first integral which is
  functionally independent together with $H$, $H_1$ and $H_2$.
\end{theorem}
Later, he developed in \cite{Ziglin:82::b} an elegant method, which 
now is called the Ziglin theory, and using it he proved
in~\cite{Ziglin:83::b} the following theorem.
\begin{theorem}[Ziglin, 1983]
\label{thm:z2}
The complexified Euler-Poisson system for a symmetric body is
integrable on $M_0$ with complex meromorphic first integrals only in
the four classical cases.
\end{theorem}
This result was improved in \cite{Ziglin:97::}.
\begin{theorem}[Ziglin, 1997]
\label{thm:z3}
The Euler-Poisson system for a symmetric body is integrable on $M_0$
with real meromorphic first integrals only in the four classical
cases.
\end{theorem}
We know that in the Euler, Lagrange and Kovalevskaya cases the
Euler-Poisson system is integrable globally, i.e. the additional first
integral exists on an arbitrary symplectic manifold $M_\chi$.  In the
Goryachev-Chaplygin case the Euler-Poisson equations are integrable
only on $M_0$, and this fact was also proved by Ziglin.

\begin{remark}
  \textsl{ We do not even try to sketch the very rich history of
    investigations of the problem of the heavy top. We refer here to
    books \cite{Golubev:53::,Leimanis:65::,Kozlov:80::,Borisov:01::}
    and references therein. However, several works where the question
    of integrability of the problem was investigated are worth
    mentioning.  In \cite{Kozlov:75::a} it was shown that for a
    non-symmetric body the Euler-Poisson equations do not admit an
    additional real analytic first integral which depends analytically
    on a small parameter $\mu$, see also \cite[Ch. III]{Kozlov:80::}.
    For a symmetric body when the ratio of the principal moments is
    small enough (i.e. for the case of the perturbed spherical pendulum)
    the non-existence of an additional real analytic first
    integral was proved \cite{Kozlov:85::c,Kozlov:86::}. A similar
    result for a perturbed Lagrange case was shown
    in~\cite{Dovbysh:90::,Bolotin:90::}. A novel, variational approach
    to proving the  non-integrability was elaborated by S.~V.~Bolotin
    in~\cite{Bolotin:92::a} where he showed the non-existence of an
    additional real analytic first integral for a symmetric heavy top
    for the case when the fixed point lies in the equatorial and the
    ratio of the principal moments of inertia is greater than 4.}
\end{remark}

The Ziglin theory is a continuation of the idea of S.~N.~Kovalevskaya
who related the (non)integrability with the behaviour of solutions of the
investigated system as functions of the complex time. The main object in
the Ziglin theory is the monodromy group of variational equations
around a particular non-equilibrium solution. As it was shown by
Ziglin, if the investigated system possesses a meromorphic first
integral, the monodromy group of variational equations possesses a
rational invariant. Thus, for an integrable system, the monodromy group
cannot be too `rich'. The main difficulty in application of the Ziglin
theory is connected with the fact that \emph{`except for a few
  differential equations, e.g., the Riemann equations,
  Jordan-Pochhammer equations and generalised hypergeometric
  equations, the monodromy group has not been determined'},
\cite[][p.~85]{Takano:76::}. For other second order 
differential equations only partial results are known, see e.g.
\cite{Churchill:90::b,Churchill:91::a}.  Having this in mind, one can
notice that the analysis of the variational equations and
determination of properties of their monodromy group given by Ziglin
in his proof of Theorem~\ref{thm:z2} is a masterpiece.  Later, the
Ziglin theory was developed and applied for a study of non-integrability of
various systems but, as far as we know, nobody used Ziglin's brilliant
technique developed in his proof of~Theorem~\ref{thm:z2}.

In the nineties the theory of Ziglin was extended by a differential
Galois approach. It was done independently by C.~Sim\'o,
J.~J.~Morales-Ruiz,
J.-P.~Ramis~\cite{Morales:94::b,Morales:99::c,Morales:01::b1,Morales:01::b2}
and A.~Braider, R.~C.~Churchill, D.~L.~Rod and M.~F.~Singer
\cite{Churchill:95::a,Churchill:96::b}. Nowadays, this approach is
called the Morales-Ramis theory. The key point in this theory is to
replace an investigation of the monodromy group of variational equations
by a study of their differential Galois group. The main fact from this
theory is similar to that of Ziglin: the existence of a meromorphic
first integral implies the existence of a rational invariant of the
differential Galois group of  variational equations. Forgetting
about differences in hypotheses in main theorems of both theories, the
biggest advantage of the Morales-Ramis theory is connected with the
fact that applying it, we have at our disposal developed tools and
algorithms of  differential algebra.  Thanks to this fact, it can be
applied more easily.

We applied the Morales-Ramis theory to study integrability of several
systems, see e.g.
\cite{Maciejewski:01::i,Maciejewski:01::j,Maciejewski:02::a,Maciejewski:02::b,Maciejewski:02::f,Maciejewski:03::a,Maciejewski:03::e},
and we notice that obtaining similar results when working only with
the monodromy group is questionable, or, at least difficult. This
gives us an idea to reanalyse the Ziglin proof of Theorem~\ref{thm:z2} which
is rather long (about 10 pages in~\cite{Ziglin:83::b}). We wanted to
present a new, much shorter and simpler proof which is based on the
Morales-Ramis theory and tools from differential algebra. In fact, at
the beginning, we believed that giving  such  proof would be a
nice and simple exercise but quickly it appeared that we were wrong. A
`naive'  application  of the Morales-Ramis theory leads quickly to
very tedious calculations or unsolvable complications. As our aim was
to give an `elementary' proof, we put a constrain on the arguments
which are allowable in it: no computer algebra. Thus we spent a lot of
time analysing sources of difficulties and complications, and the aim
of this paper is to present our own version of proofs of
Theorem~\ref{thm:z2} and Theorem~\ref{thm:z3}. As we believe, these
proofs present the whole power and  beauty of the Morales-Ramis
theory.

The plan of this paper is following. To make it self-contained in the
next section we shortly describe basic facts from the Morales-Ramis
and Ziglin theory. We collected more specific results about special
linear differential equations in Appendix. In Section~\ref{sec:pandve}
we derive the normal variational equation.  Sections~\ref{sec:Proof2}
and~\ref{sec:Proof3} contain our proofs of Theorem~\ref{thm:z2}
and~\ref{thm:z3}.  In the last section we give several remarks and
comments.
\section{Theory}
\label{sec:theory}
Below we only mention basic notions and facts concerning the Ziglin
and Morales-Ramis theory following
\cite{Ziglin:82::b,Churchill:96::b,Morales:99::c}.

Let us  consider a system of differential equations
\begin{equation}
\label{eq:ds}
\Dt x = v(x), \qquad t\in\C, \quad x\in M, 
\end{equation}
defined on a complex $n$-dimensional manifold $M$.  If $\varphi(t)$ is
a non-equilibrium solution of \eqref{eq:ds}, then the maximal analytic
continuation of $\varphi(t)$ defines a Riemann surface $\Gamma$ with
$t$ as a local coordinate.  Together with system \eqref{eq:ds} we can
also consider variational equations (VEs) restricted to $T_{\Gamma}M$,
i.e.
\begin{equation*}
\label{eq:vds}
 \dot \xi = T(v)\xi, \qquad \xi \in T_\Gamma M.
\end{equation*}
We can always reduce the order of this system by one considering the
induced system on the normal bundle $N:=T_\Gamma M/T\Gamma$ of
$\Gamma$ \cite{Ziglin:82::b}
\begin{equation}
\label{eq:gnve}
 \dot \eta = \pi_\star(T(v)\pi^{-1}\eta), \qquad \eta\in N.
\end{equation}
Here $\pi: T_\Gamma M\rightarrow N$ is the projection.  The system of
$s=n-1$ equations obtained in this way yields the so-called normal
variational equations (NVEs).  The monodromy group $\cM$ of system
\eqref{eq:gnve} is the image of the fundamental group
$\pi_1(\Gamma,t_0)$ of $\Gamma$ obtained in the process of
continuation of local solutions of \eqref{eq:gnve} defined in a
neighbourhood of $t_0$ along closed paths with the base point $t_0$.
By definition, it is obvious that $\cM\subset\mathrm{GL}(s,\C)$.  A
non-constant rational function $f(z)$ of $s$ variables $z=(z_1,\ldots,
z_s)$ is called an integral (or invariant) of the monodromy group if $
f(g\cdot z)=f(z) $ for all $g\in\cM$.

From the Ziglin theory we need the basic lemma formulated
in~\cite{Ziglin:82::b} and then given in an improved form
in~\cite{Ziglin:97::}.

\begin{lemma}
\label{lem:zig}
If system \eqref{eq:ds} possesses a meromorphic first integral defined
in a neighbourhood $U\subset M$, such that the fundamental group of
$\Gamma$ is generated by loops lying in $U$, then the monodromy group
$\cM$ of the normal variational equations  has a
rational first integral.
\end{lemma}

If system \eqref{eq:ds} is Hamiltonian, then necessarily $n=2m$ and
$M$ is a symplectic manifold equipped with a symplectic form $\omega$.
The right hand sides $v=v_H$ of \eqref{eq:ds} are generated by a
single function $H$ called the Hamiltonian of the system. For given
$H$ vector field $v_H$ is defined by $\omega(v_H,u)=\rmd H\cdot u$,
where $u$ is an arbitrary vector field on $M$. Then, of course, $H$ is
a first integral of the system. For a given particular solution
$\varphi(t)$ we fix the energy level $E=H(\varphi(t))$.  Restricting
\eqref{eq:ds} to this level, we obtain a well defined system on an
$(n-1)$ dimensional manifold with a known particular solution
$\varphi(t)$.  For this restricted system we perform the reduction of
order of variational equations.  Thus, the normal variational
equations for a Hamiltonian system with $m$ degrees of freedom have
dimension $s=2(m-1)$ and their monodromy group is a subgroup of
$\mathrm{Sp}(s,\C)$.

In the Morales-Ramis theory the differential Galois group $\cG$ of
normal variational equations plays the fundamental role, see
\cite{Morales:99::c,Morales:01::b1}.  For a precise definition of the
differential Galois group and general facts from differential algebra
see \cite{Kaplansky:76::,Ramis:90::,Beukers:92::,Magid:94::,Put:02::}.
We can consider $\cG$ as a subgroup of $\mathrm{GL}(s,\C)$ which acts
on fundamental solutions of \eqref{eq:gnve} and does not change
polynomial relations among them.  In particular, this group maps one
fundamental solution to other fundamental solutions. Moreover, it can
be shown that $\cM\subset \cG$ and $\cG$ is an algebraic subgroup of
$\mathrm{GL}(s,\C)$. Thus, it is a union of disjoint connected
components. One of them containing the identity is called the identity
component of $\cG$ and is denoted by $\cG^0$.

Morales-Ruiz and Ramis formulated a new criterion of 
non-in\-teg\-ra\-bi\-li\-ty for Hamiltonian systems in terms of the
properties of $\cG^0$ \cite{Morales:99::c,Morales:01::b1}.
\begin{theorem}
\label{thm:MR}
  Assume that a Hamiltonian system is meromorphically integrable in
  the Liouville sense in a neigbourhood of the analytic curve
  $\Gamma$. Then the identity component of the differential Galois
  group of NVEs associated with $\Gamma$   is Abelian.
\end{theorem}
In most applications the Riemann surface $\Gamma$ associated with the
particular solution is open. There are many reasons why it is better to
work with compact Riemann surfaces. Because of this, it is customary
to compactify $\Gamma$ adding to it a finite number of points at
infinity. Doing this we need a refined version of
Theorem~\ref{thm:MR}, for details see~
\cite{Morales:99::c,Morales:01::b1}. However, in the context of this
paper, the thesis of the above theorem remains unchanged if instead
of $\Gamma$ and the variational equations over $\Gamma$ we work with
its compactification.
\section{Particular solutions and variational equations}
\label{sec:pandve}
To apply the Ziglin or the Morales-Ramis theory we have to know a
non-equilibrium solution.  Let us assume that the fixed point is
located in a principal plane. Then, in fact, we can find a one
parameter family of particular solutions which describe a pendulum
like motion of the body.  For a symmetric body this assumption
is not restrictive (if necessary we can rotate the principal axes
around the symmetry axis).

 We choose the body fixed frame in the following way. Its first two
axes lie in the principal plane where the fixed point is located and
the first axis has direction from the fixed point to the centre of
mass of the body. We call this frame  special.  A map given by  
\begin{equation*}
(\bM,\bN) = (\bR \widetilde \bM, \bR \widetilde \bN),\qquad 
\bR\in\mathrm{SO}(3,\R),
\end{equation*}
transforms equations~\eqref{eq:ep} to the form 
\begin{equation*}
  \label{eq:ept}
\Dt  \widetilde \bM =[  \widetilde \bM, {\widetilde \bJ}\, {\widetilde \bM}]+ 
[ \widetilde \bN,  \widetilde \bL], \qquad
 \Dt  \widetilde \bN = [ {\widetilde \bN} ,  {\widetilde \bJ}\, {\widetilde \bM}],  
\end{equation*}
where
\[
 \widetilde \bJ =\bR^T \bJ \bR,\qquad  \widetilde \bL=\bR^T \bL.
\]
Now, if symbols with tilde correspond to the principal axes frame,
then, taking into account our assumption about the location of the
fixed point, we have
\[
 \widetilde \bJ = \begin{bmatrix} {A}^{-1} &0 &0\\
 0 & {B}^{-1} &0\\
 0 & 0& C^{-1} \end{bmatrix}, \qquad 
\widetilde \bL = [ {\widetilde L}_1, 0,{\widetilde L}_3]^T.
\]
Taking 
\[
\bR = \begin{bmatrix} 
     \phantom{-}{\widetilde L}_1 & 0 & {\widetilde L}_3\\
      0&1 &0 \\
-{\widetilde L}_3 & 0 & {\widetilde L}_1
 \end{bmatrix},
\]
we obtain 
\begin{equation*}
\bL=[1,0,0]^T, \qquad 
\bJ = \begin{bmatrix}
       a & 0& 2d\\
       0 & b & 0\\
       2d & 0 & c
\end{bmatrix},
\end{equation*}
where
\begin{gather}
\label{eq:ac}
a = \frac{{\widetilde L}_1^2}{A} +  \frac{{\widetilde L}_3^2}{C}, \qquad 
c = \frac{{\widetilde L}_1^2}{C} +  \frac{{\widetilde L}_3^2}{A}, \notag\\ 
\label{eq:d}
2 d = \left( \frac{1}{C}-\frac{1}{A} \right){\widetilde L}_1 {\widetilde L}_3, 
\qquad 
b = \frac{1}{B}.
\end{gather}
Thus, the prescribed choice of $\bR$ corresponds to the transformation
from the special to the principal frame.  Without loss of generality
we can put $b=1$. For a symmetric body we assume that $A=B\neq C$.
Under this assumption, if $d=0$, then $ {\widetilde L}_1 {\widetilde
  L}_3=0$, and, in this case, the special frame coincides with the
principal frame.

 From now on we consider the complexified Euler-Poisson system, i.e. we
assume that $(\bM,\bN)\in\C^6$.

\subsection{Case $d\neq0$}

It is easy to check that manifold
\begin{equation*}
\label{eq:inv}
\cN:= \{ ( \bM,\bN)\in \C^6\ | \ M_1=M_3=N_2 = 0, N_1^2+N_3^2=1\}\subset M_0,
\end{equation*}
is symplectic sub-manifold of $M_0$ diffeomorphic to
$T\,\bbS^1_{\C}\subset T\,\bbS^2_{\C}$ (by $\bbS^m_{\C}$ we denote
$m$-dimensional complex sphere). Moreover, $\cN$ is invariant with
respect to the flow generated by \eqref{eq:ep}.  The Euler-Poisson
equations restricted to $\cN$ have the following form
\begin{equation}
  \label{eq:rep}
  \Dt M_2 = N_3, \qquad 
  \Dt N_1=-M_2N_3, \qquad 
  \Dt N_3= M_2N_1,
\end{equation}
and are Hamiltonian with $H_{|\cN}$ as the Hamiltonian function.
For each level of Hamiltonian $H_{|\cN}=e:= 2k^2 - 1$, we obtain a 
phase curve  $\Gamma_k$.  We restrict our
attention to curves corresponding to $e\in(-1,1]$ so that $k\in(0,1]$.
A solution of system~\eqref{eq:rep} lying on the level $H_{|\cN}= 2k^2
- 1$ we denote by $(M_2(t,k),N_1(t,k),N_3(t,k))$. 

For a generic value of $k$ phase curve $\Gamma_k$ is an algebraic
curve in $\C^3\{M_2,N_1,N_3\}$, and, as intersection of two quadrics
\begin{equation*}
2k^2 - 1 = \frac{1}{2}M_2^2 + N_1, \qquad  N_1^2 + N_3^2=1, 
\end{equation*}
is an elliptic curve (for $k=1$ it is a rational curve). We can
compactify it adding two points at infinity which lie in directions
$(0,\pm\rmi,1)$. Thus, a generic $\Gamma_k$ can be considered as a
torus with two points removed. In our further consideration we work
with subfamily $\Gamma_k$ with $k\in(0,1)$. Only in
Section~\ref{sec:Proof3} we refer to the phase curve $\Gamma_1$
corresponding to $k=1$.

Equations~\eqref{eq:rep} describe the pendulum-like motions of the
body: the symmetry axis of the body remains permanently in one plane
and oscillates or rotates in it around the fixed point.  

For a point $p=(\bM,\bN)\in M_0$ by $\bv=(\bm,\bn)$ we denote a vector
in $T_p M_0$. Variational equations along phase curve $\Gamma_k$ have
the following form
\begin{equation*}
  \label{eq:varall}
\Dt
  \begin{bmatrix}
   m_1\\ m_2 \\m_3 \\ n_1 \\ n_2 \\ n_3
  \end{bmatrix}
 =
  \begin{bmatrix}
   2 d M_2 & 0 & (c-1)M_2 & 0 & 0 & 0\\ 
   0 & 0 &  0 & 0 & 0 & 1 \\
   (1-a)M_2 & 0 & -2 d M_2 & 0 & -1 & 0 \\ 
   0 & -N_3 & 0 & 0 & 0&-M_2 \\ 
 a N_3  -2 d N_1& 0 & 2 d N_3-c N_1& 0 & 0 & 0 \\ 
   0 & N_1 & 0 & M_2 & 0 & 0 
  \end{bmatrix}
  \begin{bmatrix}
   m_1\\ m_2 \\m_3 \\ n_1 \\ n_2 \\ n_3
  \end{bmatrix},
\end{equation*}
where $(M_2,N_1,N_3)=(M_2(t,k),N_1(t,k),N_3(t,k)) \in\Gamma_k$. They have
the following first integrals
\begin{equation*}
  \label{eq:ivar}
   h = M_2 m_2 + n_1 , \qquad h_1 =  N_1m_1 + N_3m_3 + M_2n_2, \qquad 
  h_2 = N_1n_1 + N_3 n_3.
\end{equation*}

As it was shown by Ziglin, the normal variational equations are given by 
\begin{equation*}
  \label{eq:nve }
\begin{split}
  \Dt m_1 &= 2 dM_2 m_1 + (c-1)M_2 m_3, \\
  \Dt m_3 & =(1-a)M_2 m_1 -2 dM_2 m_3 - n_2, \\
  \Dt n_2 &= ( a N_3  -2 d N_1) m_1 + ( 2 d N_3-c N_1 ) m_3, \\
       0 & =  N_1m_1 + N_3m_3 + M_2n_2.
\end{split}
\end{equation*}
We assume that the particular solution is not a stationary point
($M_2(t,k)\equiv0$, $N_1(t,k)\equiv\pm1$, $N_3(t,k)\equiv0$). Under this
assumption we reduce the above system to the form
\begin{equation}
  \label{eq:nve2}
\begin{split}
  \Dt m_1 &= 2 dM_2 m_1 + (c-1)M_2 m_3, \\
  \Dt m_3 & =\left[\frac{N_1}{M_2}+(1-a)M_2\right] m_1+
\left[\frac{N_3}{M_2} -2 dM_2\right] m_3.
\end{split}
\end{equation}
We can  write the above system as one second order equation 
\begin{equation}
  \label{eq:rnve}
   \Dtt m + a_1(t) \Dt m + a_0(t) m = 0, \qquad m \equiv m_1,
\end{equation}
with coefficients
\[
a_1(t) =- 2\frac{N_3(t,k)}{M_2(t,k)}, \qquad 
 a_0(t) = (1-c)N_1(t,k) + 2 d N_3(t,k) + f M_2(t,k)^2,
\]
and 
\[
f = (a-1)(c-1)-4 d^2.
\]
Now, we make the following transformation of independent variable
which is rational parametrisation of the complex circle $\bbS^1_{\C}$
\begin{equation}
\label{eq:tran}
t \longrightarrow z := \frac{N_3(t,k)}{1+N_1(t,k)}.
\end{equation}
Then we obtain
\begin{gather*}
N_1 = \frac{1-z^2}{1+z^2}, \qquad N_3 = \frac{2 z}{1+z^2},   \qquad 
M_2 = \frac{2\dot z}{1+z^2}, \\ 
{\dot z}^2 = \frac{1}{s^2+1}(z^2+1)(z^2-s^2).
\end{gather*}
where 
\begin{equation*}
s = \sqrt{\frac{1-e}{1+e}}=\frac{k'}{k}.
\end{equation*} 
After transformation~\eqref{eq:tran} equation \eqref{eq:rnve} reads  
\begin{equation}
\label{eq:var}
m''+p(z)m'+q(z)m=0,\qquad '=\dfrac{\mathrm{d}}{\mathrm{d}z},
\end{equation}
with coefficients  
\begin{equation*}
p(z)= \frac{1}{2}\left[ \frac{3}{z-\rmi} + \frac{3}{z+\rmi} -
 \frac{1}{z-s} -  \frac{1}{z+s}\right],    
\end{equation*}
 and 
\begin{equation*}
q(z) = \sum_{i=1}^4 \frac{\alpha_i}{(z-z_i)^2} + \frac{\beta_i}{z-z_i},
\end{equation*}
where we denote $(z_1, z_2,z_3,z_4)=(\rmi,-\rmi, s, -s)$,  and 
\begin{gather*}
\alpha_1 = \frac{1}{2}(1-c) - f  + \rmi d, \qquad \alpha_2 =\alpha_1^*, 
\qquad \alpha_3 =\alpha_4 = 0, \\
\beta_1 = -\frac{2d}{1+s^2} - \rmi \left(f+\frac{c-1}{1+s^2}\right),\qquad
\beta_2 = \beta_1^*,\\
\beta_3 = \frac{(1-c)(1-s^2) + 4d s}{2s(1+s^2)}, \qquad 
\beta_4 = \frac{(c-1)(1-s^2) + 4d s}{2s(1+s^2)}.
\end{gather*}
We can see that equation \eqref{eq:var} is Fuchsian and it has four
regular singular points $z_i$ over the Riemann sphere \CPOne. The
infinity is an ordinary point for this equation.  We assumed that
$k \in(0,1)$, so $s\in (0,\infty)$. Here it is important to notice
that for  real values of $a$, $c$, and $d$, and for all $s\in
(0,\infty)$ equation~\eqref{eq:var} has four distinct regular singular
points, i.e, the number of singular points does not depend on $s$. For
further calculations we fix $s=1$.

Let us note here that transformation~\eqref{eq:tran} is a branched double 
covering of Riemann sphere $\C\bbP^1\rightarrow \Gamma_k$. Moreover, the
branching points of this covering  are precisely the four points where
equation~\eqref{eq:var} has singularities.

Making the following change of the dependent  variable
\begin{equation*}
m=w\exp\left[-\dfrac{1}{2}\int_{z_0}^z p(\zeta)\mathrm{d}\zeta\right],
\end{equation*}
we can simplify \eqref{eq:var} to the standard reduced form
\begin{equation}
\label{eq:zred}
w''=r(z)w,\qquad r(z)=\dfrac{1}{2}p'(z)+\dfrac{1}{4}p(z)^2-q(z),
\end{equation}
where $r(z)$  can be
written as
\begin{equation}
\label{eq:roz}
r(z)=\sum_{i=1}^4\left[\dfrac{a_i}{(z-z_i)^2}+\dfrac{b_i}{z-z_i}\right],
\end{equation}
with coefficients
\begin{gather*}
a_2 = a_1^* =  F -\frac{1}{4}+ \rmi d , \qquad a_3=a_4= \frac{5}{16}, \\  
b_1 = b_2^* =  d +\rmi \left( F -\frac{1}{4}\right), \\
b_{3,4} =  \mp \frac{ 5}{16} -d,\qquad F = f + \frac{c}{2} - \frac{7}{16}. 
\end{gather*}
The differences of exponents $\Delta_i=\sqrt{1+4a_i}$ at 
singular point $z_i$ are the following
\begin{equation}
\label{eq:Deltai}
\Delta_1 = \Delta_2^* = \sqrt{F -\rmi d}, \qquad 
\Delta_3=\Delta_4 = \frac{3}{2}.
\end{equation}

\subsection{Case $L_3=0$}
Let us assume that $L_3=0$. Then, obviously $d=0$, and, as we already
mentioned, the special frame coincides with the principal axes frame.
As we consider a symmetric body with $A=B=1$, then additionally we
have $f=0$, and
\[
F = \frac{c}{2} -\frac{7}{16}, \qquad c =\frac{1}{C}. 
\]
Thanks to that equation~\eqref{eq:zred} has a simpler form and 
 it can be transformed to a Riemann $P$ equation.  Instead of
making a direct transformation in~\eqref{eq:zred} it is more convenient
to start from equation~\eqref{eq:rnve}. Then, instead of
transformation~\eqref{eq:tran} we make the following one
\begin{equation}
\label{eq:tran1}
t \mapsto z:=N_1(t,k)^2,
\end{equation}
and we obtain
\begin{equation}
\label{eq:rnveP}
\Dtt m + \frac{1}{2}\left( \frac{1}{z-1} + \frac{1}{2z} \right) \Dt m + 
\frac{1-c}{8} \left(\frac{1}{z-1} -\frac{1}{z} \right) m =0.
\end{equation}
For this Riemann $P$ equation the difference of exponents at $z=0$ is
$3/4$, at $z=1$ is $1/2$, and at $z=\infty$ is
\begin{equation}
\label{eq:delinf}
\Delta_\infty = \frac{1}{4}\sqrt{8c-7}. 
\end{equation}
Let us notice that as in the case $d\neq 0$
transformation~\eqref{eq:tran1} is a branched covering of Riemann
sphere $\C\bbP^1\rightarrow \Gamma_k$, and the branching points of
this covering are precisely the three points where
equation~\eqref{eq:rnveP} has singularities.
\subsection{Case  $L_3=0$. Second particular solution.}
When $L_3=0$ we have at our disposal another family of particular solutions.
Under our assumption $A=B=1$ and $\bL=[1,0,0]^T$, the following manifold
\begin{equation*}
\cN_1:= \{ ( \bM,\bN)\in \C^6\ | \ M_1=M_2=N_3 = 0, N_1^2+N_2^2=1\}\subset M_0,
\end{equation*}
is invariant with respect to the flow of system~\eqref{eq:ep}.
Similarly as $\cN$, manifold $\cN_1$ is diffeomorphic to
$T\,\bbS^1_{\C}\subset T\,\bbS^2_{\C}$ and it is a symplectic
sub-manifold of $M_0$. The Euler-Poisson equations restricted to
$\cN_1$, read
\begin{equation}
  \label{eq:rep1}
  \Dt M_3 =- N_2, \qquad 
  \Dt N_1= c M_3N_2, \qquad 
  \Dt N_2=-c M_3N_1.
\end{equation}
We consider a family of phase curves $k\mapsto \Gamma^1_k$ of the above
equations given by
\begin{equation}
\label{eq:lev1}
    \frac{1}{2}cM_3^2 + N_1 = e, \qquad N_1^2 + N_2^2 = 1,
\end{equation}
where $e=2k^2-1$.  For $k\in (0,1)$ curves $\Gamma_k^1$ are
non-degenerate elliptic curves.  Variational equations along phase
curve $\Gamma_k^1$ have the following form
\begin{equation*}
  \label{eq:var1}
\Dt
  \begin{bmatrix}
   m_1\\ m_2 \\m_3 \\ n_1 \\ n_2 \\ n_3
  \end{bmatrix}
 =
  \begin{bmatrix}
   0 &  (c-1)M_3 &0 & 0 & 0 & 0\\ 
   (1-c)M_3 & 0 &  0 & 0 & 0 & 1 \\
   0 & 0 & 0 & 0 & -1 & 0 \\ 
   0 & 0 & cN_2 & 0 & cM_3 & 0 \\ 
   0 & 0 & -cN_1& -cM_3& 0 & 0 \\ 
   -N_2 & N_1 & 0 & 0 & 0 & 0 
  \end{bmatrix}
  \begin{bmatrix}
   m_1\\ m_2 \\m_3 \\ n_1 \\ n_2 \\ n_3
  \end{bmatrix},
\end{equation*}
and they have the following first integrals
\begin{equation*}
  \label{eq:ivar1}
\begin{split}
   h &= cM_3 m_3 + n_1 , \\
  h_1& =  N_1m_1 + N_2m_2 + M_3n_3, \\
  h_2 &= N_1n_1 + N_2 n_2.
\end{split}
\end{equation*}
The normal variational equations are given by 
\begin{equation*}
\label{eq:nve1}
\begin{split}
\Dt m_1 &= (c-1)M_3 m_2, \\
\Dt m_2 &= (1-c)M_3 m_1 + n_3, \\
\Dt n_3 &= -N_2 m_1 +N_1 m_2 , \\
 0&=  N_1m_1 + N_2m_2 + M_3n_3.
\end{split}
\end{equation*}
Now, the reduction of the above system to the second order equation gives 
\begin{equation}
\label{eq:nves1}
\ddot n + (M_3^2 - N_1) n =0, \qquad n\equiv n_3.
\end{equation}
Let us notice that from equations~\eqref{eq:rep1} and \eqref{eq:lev1}
it follows that
\begin{equation*}
\label{eq:dN1}
{\dot N_1}^2 = 2c(e-N_1)(1-N_1^2).
\end{equation*} 
Thus, putting 
\[
N_1=\frac{2}{c} v + \frac{e}{3},
\]
we obtain the following equation
\begin{equation*}
\label{eq:vwp}
{\dot v}^2 = 4v^3 -g_2 v - g_3,
\end{equation*}
determining the Weierstrass function $v(t)=\wp(t;g_2,g_3)$ with invariants 
\begin{equation*}
\label{eq:g2g3}
g_2 = \frac{1}{3}c^2(e^2+3), \qquad  g_3 = \frac{1}{27}c^3 e (e^2 -9) .
\end{equation*}
Hence, we can express $N_1$, and $M_3^2$ (using \eqref{eq:lev1}), in
terms of the  Weierstrass function $\wp(t;g_2,g_3)$. The discriminant and
the modular  function of $\wp(t;g_2,g_3)$ are following
\begin{equation*}
\label{eq:dedi}
\begin{split}
&\Delta = g_2^3 -27g_3^2 = c^6(e^2-1)^2, \\
&j(g_2,g_3) = j(e)=\frac{g_2^3}{g_2^3 - 27g_3^2}=\frac{(e^2+3)^3}{27(e^2-1)^2}.
\end{split}
\end{equation*}
Hence, we can rewrite equation~\eqref{eq:nves1} in the form
of the Lam\'e equation
\begin{equation}
\label{eq:nve1lam}
   \Dtt n = (\alpha \wp(t;g_2,g_3) + \beta) n , 
\end{equation}
where
\[
\alpha = 2C(2C+1), \qquad \beta = \frac{e}{3}C(1-4C).
\] 
It is important to notice here the physical restriction on parameter
$C$, namely,  we have $C\in(0,2)$.

\section{Proof of Theorem~\ref{thm:z2}}
\label{sec:Proof2}
In our proof of Theorem~\ref{thm:z2} we try to be as close as possible
to the proof of Ziglin. Namely, first we show that a necessary condition
for integrability is $\widetilde L_3=0$ (or $\widetilde L_1=0$, but
this gives the already known integrable case of Lagrange). In fact this is
the most difficult part of the proof. Then, we use the second family
of particular solutions and we restrict the possible values of the
principal moment of inertia. Finally, using the first solution, we
limit all allowable values of $C$ to those corresponding to the known
integrable cases.

We organise the three steps of the proof  in the form of three lemmas.
Only the first one is somewhat involved, the remaining two are very
simple.

The first step is to show that a necessary condition for 
integrability is $d=0$, see formulae~\eqref{eq:d}.

\begin{lemma}
\label{lem:1}
Let us assume that $d\neq 0$. Then the identity component of the
differential Galois group of equation~\eqref{eq:zred} is not Abelian.
\end{lemma}
\begin{proof}
  In our proof we  use of Lemma~\ref{lem:alg} and~\ref{lem:algc1},
  see Appendix. If~\eqref{eq:zred} is reducible then the identity
  component $\cG^0$ of its differential Galois group is Abelian in two
  cases: when $\cG$ is a subgroup of diagonal group $\cD$ or when it
  is a proper subgroup of triangular group $\cT$.  
  
  First we show that $\cG\not\subset\cD$. Let us assume the opposite. Then
  there exist two exponential solutions of~\eqref{eq:zred} which have
  the following form
\begin{equation*}
\label{eq:expsol}
 w_l = P_l \prod_{i=1}^4(z-z_i)^{e_{i,l}}, \qquad P_l\in\C[z], \quad l=1,2, 
\end{equation*}
where $e_{i,l}$ for $l=1,2$ are exponents at singular point $z_i$, i.e., 
\[
e_{i,l}\in\left\{ \frac{1}{2}(1 + \Delta_i), \frac{1}{2}(1 -\Delta_i)\right\}.
\]
Here $\Delta_i$ for $i=1,\ldots, 4$ are given by \eqref{eq:Deltai}. 
The product of these solutions $v=w_1w_2$ belongs to $\C(z)$ and it is
a solution of the second symmetric power of~\eqref{eq:zred}, i.e.
equation~\eqref{eq:ssp} with $r$ given by~\eqref{eq:roz}. This
equation has the same singular points as equation~\eqref{eq:zred}.
Exponents $\rho_{i,l}$ at singular points $z_i$, and at infinity
$\rho_{\infty,l}$ for the second symmetric power of~\eqref{eq:zred}
are given by
\[
\rho_{i,l}\in  \{1, 1\pm\Delta_i\}, \qquad 
\rho_{\infty,l}\in\{-2,-1,0\}, \qquad l=1,2,3, 
\]
where $\Delta_i$ for $i=1,2,3,4$ are given by~\eqref{eq:Deltai}.  
If we write $v = P/Q$ with $P,Q \in \C[z]$ then
\[
Q = \prod_{i=1}^K(z-r_i)^{n_i}, \qquad n_i\in\N, \qquad
r_i\in\{z_1,z_2,z_3,z_4\},
\]
and $ n_i= -\rho_{i,l}\in\N$ for certain $l$. However, if $d\neq0$,
then $\rho_{i,l}$ is not a negative integer for $i=1,2,3,4$ and
$l=1,2,3$. This implies that $Q=1$. Hence, equation~\eqref{eq:ssp} has
a polynomial solution $v=P$, and $\deg P = -\rho_{\infty,l}\leq 2$.
But $v$ is a product of two exponential solutions of the form
~\eqref{eq:expsol}, so we also have
\begin{equation}
\label{eq:exprod}
v = P_1P_2 \prod_{i=1}^4(z-z_i)^{e_{i,m}+e_{i,l}}\in\C[z], 
\qquad m,l\in\{1,2\}.
\end{equation}
Consequently, ${e_{i,m}+e_{i,l}}$ is a non-negative integer for
$i=1,2,3,4$.  As for $d\neq0$, we have $2e_{i,l}\not\in\Z$ for
$i=1,2,3,4$ and $l=1,2$, we deduce that in~\eqref{eq:exprod} $m\neq
l$.  But $e_{i,1}+e_{i,2}=1$, for $i=1,2,3,4$.
Thus, we have 
\[ 
\deg v = \deg P =
\deg(P_1P_2)+4\geq 4.
\]
We have a contradiction because we already showed  that $\deg P \leq 2$.

It is also impossible that $\cG$ conjugates to $\cT_m$ for a certain
$m\in\N$ because when $d\neq 0$ exponents for $z_1$ and $z_2$ are not
rational. This implies also that $\cG$ is not finite.

The last possibility that $\cG^0$ is Abelian occurs when $\cG$ is
conjugated with a subgroup of $D^\dag$. We show that it is impossible.
To this end, we apply the second case of the Kovacic algorithm, see
Appendix. The auxiliary sets for singular points are following
\[
E_1=E_2=\{2\}, \qquad E_3=E_4 = \{-1,2,5\}, \qquad E_\infty=\{0,2,4\}.
\]
In the Cartesian product $E=E_\infty\times E_1 \times \cdots \times E_4$ 
we look for such  elements $e$ for which 
\[
d(e):=\frac{1}{2}\left( e_\infty - \sum_{i=1}^4 e_i\right),
\]
is a non-negative integer. There are two such elements, namely
\[
e^{(1)}=(4,2,2,-1,-1),\qquad e^{(2)}=(2,2,2,-1,-1).
\]
We have $d(e^{(1)})=1$ and $d(e^{(2)})=0$. We have to check if there
exists polynomial $P = p_1 z + p_0$ which satisfies the following
equation
\[ 
P''' + 3\theta P'' +(3 \theta^2 + 3\theta' -4r)P' + 
(\theta'' + 3 \theta\theta' + \theta^3 -4 r\theta - 2r')P =0,
\] 
where
\[
\theta = \frac{1}{z-\rmi} +  \frac{1}{z+\rmi} -\frac{1}{2}\frac{1}{z-1}-
\frac{1}{2}\frac{1}{z+1}. 
\]
Inserting $P$ into the above equation we obtain the following system of
linear equations for its coefficients
\begin{gather*}
- d p_0 + F p_1 = 0 , \qquad (1-2F)p_0 + 6 d p_1 =0 \\
F p_0 + d p_1 = 0 , \qquad -6 d p_0 + (1-2F) p_1 =0.
\end{gather*}
The above system for $p_0$ and $p_1$ has a non-zero solution if $d^2 +
F^2=0$, but for a real $d$ and $F$ it is possible only when $d=0$ and
$F=0$.
\end{proof}
As the covering $t\mapsto z$ given by~\eqref{eq:tran} does not change
the identity component of the differential Galois group of the normal
variational equations~\eqref{eq:nve2}, from the above lemma it follows
that if the Euler-Poisson equations are integrable, then
\[ 
2 d = \left( \frac{1}{C}-\frac{1}{A} \right){\widetilde L}_1 {\widetilde L}_3 =0.
\]
For a symmetric body when $A=B\neq C$, the above condition implies
that either ${\widetilde L}_1 = 0$, and this corresponds to the
integrable case of Lagrange, or ${\widetilde L}_3=0$. Hence, we have
to investigate the last case.  Notice that in this case the special
frame is the principal axes frame so we have ${\widetilde
  \bL}=\bL=[1,0,0]^T$, $a=1/A=1/B=b=1$ and $c=1/C\neq 1$.  At this
point it is worth to observe that now the identity component of the
differential Galois group of equation~\eqref{eq:zred} is Abelian for
infinitely many values of $c$. Let us remind  here that from the
physical restriction it follows only that $c\in (1/2,\infty )$.
\begin{proposition}
\label{prop:1}
If $L_3=0$ then the identity component of the
differential Galois group of equation~\eqref{eq:zred} is Abelian in the 
following cases:
\begin{gather}
c = 1 + 2l(4l-1), \quad  c= 4 + 2l(4l+5), \quad
 c = 2 +2l(4l+3,)\\
 c= \frac{11}{8} +2l(l+1), \qquad c = \frac{79}{72} + \frac{4}{3}l + 2l^2,
\end{gather}
where $l$ is an integer. 
\end{proposition}
\begin{proof}
  When $L_3=0$ then the identity component of the differential Galois
  group of equation~\eqref{eq:zred} is Abelian if and only if the
  differential Galois group of equation~\eqref{eq:rnveP} is Abelian
  (transformation from \eqref{eq:zred} to \eqref{eq:rnveP} is
  algebraic). Then applying Kimura Theorem~\ref{thm:kimura} to
  equation~\eqref{eq:rnveP} we easily derive the above values of $c$  for
  which the identity component of the differential Galois
  group of this equation is Abelian.
\end{proof}
Thus, applying the Morales-Ramis or Ziglin theory and using the first
particular solution we cannot prove non-integrability of the
Euler-Poison equations for all values of $c$ listed in the above
proposition.  This is why, in the lemma below, we consider normal
variational equations corresponding to the second particular solution.

\begin{lemma}
\label{lem:2}
  Assume that  $C\in(0,2)$ and $C\neq m/4$, for $m=1,\ldots, 7$, then
  for almost all $e\in\R$, the differential Galois group of
  equation~\eqref{eq:nve1lam} is \SLtwoC.
\end{lemma}  
\begin{proof}
  We assume first that $e\neq\pm 1$. Then the discriminant of the elliptic
  curve associated with $\wp(t;g_2,g_3)$ does not vanish, and we can
  apply Lemma~\ref{lem:lame}, see Appendix.  We consider successively
  three cases from this lemma.
  
  For the Lam\'e-Hermite case,  we have $\alpha = n(n+1)$ for $n\in\Z$.
  This implies that $C=n/2$, and hence, as $C\in(0,2)$ we have $C\in\{1/2,1,3/2\}$. 
  
  For the Brioschi-Halphen-Crawford case,  we have $\alpha = n(n+1)$ and
  $m=n+1/2 \in \N$. Thus, we have
\[
  C = -\frac{1}{4} +\frac{1}{2}m, \qquad m\in\N.
\]
So, this case can occur only when $C\in\{1/4,3/4,5/4,7/4\}$. 

In the Baldassarri case we notice that the mapping
\[
 \R\backslash\{-1,1\} \ni e \mapsto  j(e),
\]
is non-constant and continuous. Hence, by Dwork
Proposition~\ref{prop:dwork}, for a fixed $C$ this case can occur only
for a finite number of values of $e$.
\end{proof} 

As $C=1/2$, $C=1$ and $C=1/4$ correspond to the Kovalevskaya, Euler
and Goryachev-Chaplygin cases, respectively, we have to investigate
cases 
\[
C\in\{3/4,5/4,3/2,7/4\}.
\]
To this end we return to equation~\eqref{eq:zred}. As we show, for
$d=0$ it can be transformed to the form~\eqref{eq:rnveP} and, moreover,
this transformation does not change the identity component of its
differential Galois group. We can prove the following.
\begin{lemma}
\label{lem:3}
  For $C\in\{3/4,5/4,3/2,7/4\}$ the  differential Galois group
  of~\eqref{eq:rnveP} is \SLtwoC.
\end{lemma}
\begin{proof}
  For $C\in\{3/4,5/4,3/2,7/4\}$ the respective values of the difference of
  exponents at infinity $\Delta_\infty$ (see
  formula~\eqref{eq:delinf}) for equation~\eqref{eq:rnveP} are
  following
\[
\frac{1}{4}\left\{ \sqrt{\frac{11}{3}}, \rmi  \sqrt{\frac{3}{5}},
\rmi  \sqrt{\frac{5}{3}}, \rmi  \sqrt{\frac{17}{7}}\right\}.
\] 
Now, a direct inspection of possibilities in the Kimura
Theorem~\ref{thm:kimura} shows that Riemann $P$
equation~\eqref{eq:rnveP} with prescribed differences of exponents does
not possess a Liouvillian solution, so its differential Galois group
is \SLtwoC.
\end{proof}
Now the proof of Theorem~\eqref{thm:z2} is a simple consequence of the
above three lemmas. 
\section{Proof of Theorem~\ref{thm:z3}}
\label{sec:Proof3}
In the proof of Theorem~\ref{thm:z3} we apply the Ziglin
Lemma~\ref{lem:zig} and his idea of its application given in
\cite{Ziglin:97::}. 

The Euler-Poisson equations restricted to $\cN$ possess a
hyperbolic equilibrium at $u=(0,0,0,1,0,0)$. The phase curve
$\Gamma_1$ corresponds to the solution of equations~\eqref{eq:rep} with
$k=1$. It contains two real components which are real phase curves
corresponding to real solutions homoclinic to $u$. Their union is $\Re
\Gamma_1$ and we denote its closure by $\Omega$. 
\begin{lemma}
\label{lem:eps}
  For an arbitrary complex neighbourhood $U\subset\cN$ of $\Omega$ there
  exists $\epsilon>0$, such that for $0<1-k<\epsilon$ the fundamental
  group $\pi_1(\Gamma_k)$ of phase curve $\Gamma_k$ is generated by
  loops lying in $U$.
\end{lemma}
\begin{proof}
The time parametrisation
of $\Gamma_k$ is given by
\begin{equation}
\label{eq:exfik}
\begin{split}
M_2(t,k)&:=-2k\cn( t,k), \\
N_1(t,k)&:= 2k^2\sn^2( t,k)-1,\\
N_3(t,k)&:= 2k\sn( t,k)\dn( t,k), 
\end{split}
\end{equation}
where $\sn(t,k)$, $\cn(t,k)$ and $\dn(t,k)$ denote the
Jacobi elliptic functions of argument $t$ and modulus $k$.
Thus,  particular solutions of~\eqref{eq:ep}
\begin{equation}
\label{eq:phitk}
\varphi(t,k):= (0,M_2(t,k),0,N_1(t,k),0,N_3(t,k)),
\end{equation}
defined by~\eqref{eq:exfik} are single-valued, meromorphic, and double
periodic with periods
\begin{equation}
\label{eq:Tik}
  T_1(k)= 2K(k)+2\rmi K'(k),\qquad T_2(k)= 2K(k)-2\rmi K'(k),
\end{equation}
where $K(k)$ is the complete elliptic integral of the first kind with
modulus $k$, $K'(k):=K(k')$, and $k':=\sqrt{1-k^2}$. In each period
cell they have  two simple poles at:
\begin{equation*}
 t_1(k)= \rmi K'(k), \quad t_2(k)= -\rmi K'(k) \mod (T_1(k),T_2(k)).
\end{equation*}
Periods $T_1(k)$ and $T_2(k)$ given by~\eqref{eq:Tik} of
  solution~\eqref{eq:phitk} are primitive.  Minimal real and imaginary
  periods are $T(k)=4K(k)$ and $T'(k)=4\rmi K'(k)$.  As a base point
  $x(k)\in\Gamma_k$ we choose $x(k)=\varphi(t_0(k),k)$ where
  $t_0(k)=K(k)$. Let us notice that from \eqref{eq:exfik} it follows
  that
\begin{equation}
\label{eq:xk}
 \begin{split}
 M_2(t_0(k),k) &=0, \\
 N_1(t_0(k),k)&=2k^2-1, \\ 
 N_3(t_0(k),k)&=2kk'. 
\end{split}
\end{equation}
Now, let 
\[
 \lambda_k, \lambda_k':[0,1]\rightarrow\Gamma_k,
\]
be the loops with base point $x(k)$ corresponding to periods $T(k)$
and $T'(k)$, respectively. These loops cross at point 
\begin{equation*}
\label{eq:xkp}
x'(k)=\varphi(t_0(k)+T(k)/2,k)= \varphi(t_0(k)+T'(k)/2,k ). 
\end{equation*}
As a results, we obtain four semi-loops with end points $x(k)$ and
$x'(k)$. The fundamental group $\pi_1(\Gamma_k,x(k))$ of $\Gamma_k$ is
generated by these semi-loops, see Figure~\ref{fig:loops}.
\begin{figure}[h]
\centering \includegraphics[scale=0.8]{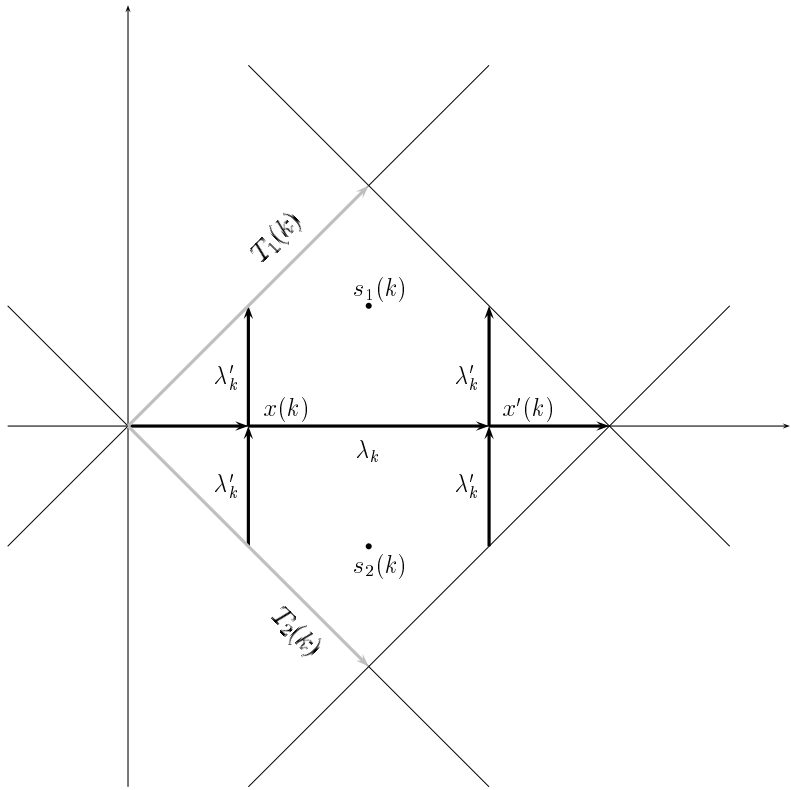}
\caption{\label{fig:loops}Parallelogram of period with marked loops}
\end{figure}
Let us analyse what happens when $k$ tends to 1.  From \eqref{eq:xk}
it follows that $x(k)$ tends to $u$ and from \eqref{eq:exfik} we
deduce that loop $\lambda_k$ tends to $\Omega$.  To see what happens
with loop $\lambda'_k$ when $k$ tends to 1, let us put $t=t_0(k)+\rmi
\tau$ in formulae \eqref{eq:exfik}. We obtain
\begin{equation*}
\begin{split}
M_2(t,k)&= 2\rmi kk'\frac{\sn(\tau,k')}{\dn(\tau,k')},\\
N_1(t,k)&= -1+\frac{2k^2}{\dn^{2}(\tau,k')} \\
N_3(t,k)&= 2kk'\frac{\cn(\tau,k')}{\dn^{2}(\tau,k')}.
\end{split}
\end{equation*}
Thus, loop $\lambda'_k$ tends to point $u$ as $k$ tends to 1. 
\end{proof}

Time parametrisation of these
phase curves $\Gamma_k^1$ is given by
\begin{equation}
\label{eq:sol2t}
\begin{split}
M_3(t,k)&=\frac{2k}{\omega}\cn(\omega t,k),\\
N_1(t,k)&=2k^2\sn^2(\omega t,k)-1, \\
N_2(t,k)&=2k\sn(\omega t,k)\dn(\omega t,k),
\end{split}
\end{equation}
where $\omega^2 = c$. Thus, the second family of particular solutions
of~\eqref{eq:ep}
\begin{equation*}
\label{eq:phi1tk}
\varphi_1(t,k):= (0,0,M_3(t,k),N_1(t,k),N_2(t,k), 0),
\end{equation*}
defined by~\eqref{eq:sol2t} contains solutions which are
single-valued, meromorphic, and double periodic with periods
\begin{equation*}
\label{eq:Tik1}
  T_1(k)= 2\frac{K(k)}{\omega}+2\rmi \frac{K'(k)}{\omega},\qquad 
T_2(k)= 2\frac{K(k)}{\omega}-2\rmi \frac{K'(k)}{\omega}.
\end{equation*}
In each period cell they have two simple poles at:
\begin{equation*}
 t_1(k)= \rmi \frac{K'(k)}{\omega}, \quad 
t_2(k)= -\rmi \frac{K'(k)}{\omega} \mod (T_1(k),T_2(k)).
\end{equation*}
Let us notice that for the Euler-Poisson equations restricted to
$\cN_1$ point $u=(0,0,0,1,0,0)$ is also a hyperbolic
equilibrium.  The phase curve $\Gamma_1^1$ corresponds to the solution of
equations~\eqref{eq:rep1} with $k=1$. As in the previous case, it
contains two real components which are real phase curves corresponding
to real solutions homoclinic to $u$. Their union is $\Re \Gamma_1^1$
and we denote its closure by $\Omega_1$. Using the same arguments as
in the proof of Lemma~\ref{lem:eps} we can show the following.
\begin{lemma}
\label{lem:eps1}
  For an arbitrary complex neighbourhood $U\subset\cN_1$ of $\Omega_1$ there
  exists $\epsilon>0$, such that for $0<1-k<\epsilon$ the fundamental
  group $\pi_1(\Gamma_k^1)$ of phase curve $\Gamma_k$ is generated by
  loops lying in $U$.
\end{lemma}
Now, to prove Theorem~\ref{thm:z3} let us notice that we showed that,
except for the known integrable cases, the identity component of the
differential Galois group of the normal variational equations
corresponding to $\Gamma_k$ or $\Gamma^1_k$ is not Abelian for almost
all values of $k\in(0,1)$.  In fact, in Lemma~\ref{lem:1} we proved
that for $d\neq 0$ the identity component of the differential Galois
group of the normal variation equations corresponding to $\Gamma_k$
with $k=1/\sqrt{2}$ is not Abelian. By Lemma~\ref{lem:par} it is not
Abelian for almost all values of $k\in (0,1)$.  Then,
in~Lemma~\ref{lem:3} we proved that the identity component of the
differential Galois group of the normal variation equations
corresponding to $\Gamma_k^1$ is not Abelian for almost all values of
$k\in (0,1)$, except for $C=m/4$, $m=1,\ldots, 7$.  Finally, for
$C=m/4$ such that $C\not\in\{1/4, 1/2,1\}$, we showed that the identity
component of the differential Galois group of the normal variation
equations corresponding to $\Gamma_k$ with $k=1/\sqrt{2}$ is not
Abelian. Again, by Lemma~\ref{lem:par} it is not Abelian for almost
all values of $k\in (0,1)$.  

Both normal variational equations corresponding to $\Gamma_k$ and
$\Gamma_k^1$ are Fuchsian.  For a Fuchsian equation we know that if
the identity component of its differential Galois group is not
Abelian then its monodromy group does not possess a rational
invariant, see Theorem~3.17 in \cite{Churchill:96::b}.

Assume now that for $C\not\in\{1/4,1/2,1\}$ the Euler-Poisson
equations possess an additional real meromorphic first integral
defined in a real neighbourhood of $ \Omega\cup\Omega_1$. Then we
can extend this integral to a complex meromorphic one, defined in a
certain complex neighbourhood $U$ of $ \Omega\cup\Omega_1$.  Then, by
Lemma~\ref{lem:eps} and \ref{lem:eps1}, we find such $\epsilon>0$
that the fundamental groups of $\Gamma_k$ and $\Gamma^1_k$ with
$0<1-k<\epsilon$, are generated by loops lying entirely in $U$. Then,
from the Ziglin Lemma~\ref{lem:zig}, it follows that both monodromy groups
of normal variational equations corresponding to $\Gamma_k$ and
$\Gamma^1_k$ possess a rational invariant. However, above we showed that
at least for one of them it is not true. A contradiction proves
Theorem~\ref{thm:z3}. 

\section{Remarks and Comments} 

One important difference between the Ziglin and Morales-Ramis theory
is related with the procedure of obtaining the normal variational
equations.  Assume that system~\eqref{eq:ds} possesses certain number
of known first integrals $H_i$, such that their differentials $\rmd
H_i$, are linearly independent on $\Gamma$. Then $\rmd H_i\circ
\pi^{-1}$ for $i=1,\ldots,k$ are independent first integrals of
\eqref{eq:gnve}.Their common level
\begin{equation*}
N_p:= \{ \eta\in F\,|\,   \rmd H_i\circ \pi^{-1}\eta = p_i,\quad p_i\in \C,
 \quad i=1,\ldots,k\},
\end{equation*}
defines a $m$-dimensional linear bundle over $\Gamma$, where
$m=n-k-1$.  Using these integrals we can reduce the order of system
\eqref{eq:gnve}. Namely, we consider the reduced normal variational
equations
\begin{equation}
\label{eq:rgnve}
 \dot \eta = \pi_\star(T(v)\pi^{-1}\eta), \qquad \eta\in N_p.
\end{equation}
However after this reduction defined by Ziglin, instead of a linear,
we have an affine bundle over $\Gamma$, equations \eqref{eq:rgnve} are
generally not homogeneous ones, and the monodromy group is a subgroup
of affine transformations of $\C^m$.  Till now this construction has
not been translated to the Morales-Ramis theory where we work with a
system of homogeneous equations defined on $N_0$. To realise the
importance of Ziglin reduction let us notice that using only the
reduced equation on $N_0$ it is impossible to prove global
non-integrability of the Goryachev-Chaplygin case. In
\cite{Ziglin:97::} he gave such a non-integrability proof
investigating the reduced normal variational equations on $N_p$ with
$p\neq 0$.

The Morales-Ramis theory is coordinate independent, however,
investigating a specific problem, we always have to choose appropriate
coordinates. The form of normal variational equations depends on local
coordinates and this is why their choice is important. It is
especially evident when we investigate a problem connected with a
rigid body. Equations of motion of the heavy top can be written in
many different forms.  As we mentioned in Remark~\ref{rem:r1}, the
natural phase space for a rigid body with a fixed point is
$T^*\mathrm{SO}(3,\R)$. There are no `natural' coordinates on
$\mathrm{SO}(3,\R)$, and thus there are no `natural' canonical
coordinates. The most widely used are Androyer-Deprit canonical
coordinates~\cite{Borisov:01::} or the Euler angles and conjugated
momenta. In fact we checked which, from almost all known coordinates
on $T^*\mathrm{SO}(3,\R)$, are most feasible for application of the
Morales-Ramis theory.   In our exposition we
work with the Euler-Poisson equations.  However, our choice of the
body fixed frame is not conventional.  Usually the principal axes
frame is used. To see what is an advantage of our choice, let us
notice that using the principal axes frame we can derive the normal
variational equation in the form similar to~\eqref{eq:rnve}, however,
to put it in the form of an equation with rational coefficient we have
to choose a transformation different than~\eqref{eq:tran}, and as a
result, we obtain, instead of equation~\eqref{eq:zred} possessing four
regular singularities, a much more complicated Fuchsian equation with
seven singular points. Our choice of the body fixed frame appears e.g.
in \cite{Dokshevich:92::}.

Simplifications of the normal variational equations which occur when
$L_3=0$ need an explanation. In fact, one can observe that although
the Riemann surface $\Gamma_k^1$ for the second particular solution is
a torus with two points removed, see formulae~\eqref{eq:sol2t}, the
normal variational equation~\eqref{eq:nve1lam} corresponding to it
has the form of a Lam\'e equation, so it is defined over a torus with
one point removed.  The reason of what happened is  symmetry. When
$L_3=0$ the Euler-Poisson equations restricted to $M_0$ are invariant
with respect to an involutive symplectic diffeomorphism $\cJ_1:
M_0\mapsto M_0$ defined by
\[
  \cJ_1(M_1,M_2,M_3, N_1,N_2,N_3) = (-M_1,M_2,-M_3, N_1,-N_2,N_3). 
\]
Let us denote
\[
M=\{ p\in M_0\,|\, \cJ_1(p)\neq p\}, \qquad \widehat M=
M/\cJ_1,
\]
and let $\pi:M\mapsto \widehat M$ be the projection. In the natural way
equations~\eqref{eq:ep} induce Hamiltonian equations on $\widehat M$
with Hamiltonian function $\widehat H= H\circ\pi^{-1}$. Then,
according to Ziglin, see Lemma on page 36 in \cite{Ziglin:82::b}, if
system~\eqref{eq:ep} is integrable, then the induced Hamiltonian system
on $\widehat M$ is also integrable. For the induced system we have
a family of particular solutions $\widehat\varphi_1(t,k) =
\pi\circ\varphi_1(t,k)$. The corresponding Riemann surfaces
$\widehat\Gamma^1_k$ are tori with one point removed.

A simplification of the normal variational equation for solution
$\varphi(t,k)$ when $L_3=0$ and the fact that we can transform them to a
Riemann $P$ equation (for an appropriate choice of  energy) is also
related with symmetry. Namely, when $L_3=0$, system \eqref{eq:ep}
restricted to $M_0$ is also invariant with respect to an involutive
symplectic diffeomorphism $\cJ: M_0\mapsto M_0$ defined by
\[
  \cJ(M_1,M_2,M_3, N_1,N_2,N_3) = (-M_1,-M_2,M_3, N_1,N_2,-N_3). 
\]
For symmetry reduction of variational equations see Section~4.2 in
\cite{Churchill:96::b}.

Let us note that in Lemma~\ref{lem:1} we claim that if $d\neq 0$ then
the identity component of the differential Galois group of
equation~\eqref{eq:zred} is not Abelian. Thus, it can be the whole
group \SLtwoC\ or whole triangular subgroup $\cT$ of \SLtwoC. We do not
know if the second case can occur.

In the case of first particular solution we do not work with the
elliptic curve $\Gamma_k$ but with the Riemann sphere (minus singular
point) for which $\Gamma_k$ is a covering. The reason of this is that
we have no tool similar to the Kovacic algorithm for a second order
linear differential equation defined on an elliptic curve. However, in
the case of the second particular solution we can work directly on
elliptic curve $\Gamma^1_k$ because, in this case, the normal
variational equation is the Lam\'e equation for which the monodromy
group is know. Of course, in this case we can also work on the Riemann
sphere making well know transformation of the Lam\'e equation to its
algebraic form.

In his proof of Theorem~\ref{thm:z2} and \ref{thm:z3} Ziglin used the
explicit time parametrisation of particular solutions. First he showed
that if the system is integrable then the monodromy of the normal
variational equations along real periods of a particular solution must
be equal to the identity. Then, using analytical tools he derived the
necessary conditions for the integrability. 

In our exposition we use the explicit time parametrisation of
particular solutions in the proof of Theorem~\ref{thm:z3}. In fact, we
use it only to show explicitly what happens with loops along real and
imaginary periods when $k$ tends to 1. However, one can deduce this
information from the equations defining the elliptic curve. Thus, we
can avoid using explicit time parametrisation at all. We keep it in
the proof of Theorem~\ref{thm:z3} because, as we hope, it makes the
exposition more transparent. 

The physical restriction $C\in(0,2)$ plays crucial role in our, as
well as, in the  Ziglin proof. Integrable systems are really rare, hence it
is an interesting question if the Euler-Poisson equations are
integrable for values of parameters which do not satisfy this restriction.

Considering the case $L_3=0$ and the first particular solution Ziglin
showed that the necessary condition for integrability is $c\in\N$. In
our Proposition~\ref{prop:1} there are two families of $c$ such that
$c\in\Q$. The reason why they appears is that we fixed the energy for
the first solution.

In the proof of Lemma~\ref{lem:2} we show
that the Brioschi-Halphen-Crawford case is possible only when
$C\in\{1/4,3/4,5/4,7/4\}$. For this values of $C$ we can calculate the
Brioschi determinant $Q_m(g_2,g_3,\beta)$ defined by~\eqref{eq:bdet}.
Calculations show that it vanishes identically only when $C=1/4$,
i.e., for the Goryachev-Chaplygin case. Thus, in fact, to prove that
for $C\in\{3/4,5/4,7/4\}$ the Euler-Poisson equations are
non-integrable, we can use the second solution. We use the first one 
because calculations are simpler.

\section{Acknowledgements}
We are very thankful to Mich\`ele Audin for her remarks, comments,
suggestions and corrections. They allowed us to improve considerably
not only the contents of the paper, but also gave a more clear proof
of our main result.

We would like to thank Delphine Boucher, Juan J.~Morales-Ruiz,
Jacques-Arthur Weil, Carles Sim\'o, Michael F.~Singer and Felix Ulmer
for discussions and help which allowed us to understand many topics
related to this work. We are very grateful to Mark van Hoeij with whom
we started to discuss some problems concerning this paper in the end
of the previous century. Many important comments by Robert~S.~Maier
concerning the Lam\'e equation are gratefully acknowledged.  

As usual, we thank Zbroja (Urszula Maciejewska) not only for her
linguistic help.  

For the second author this research has been supported by a Marie
Curie Fellowship of the European Community programme Human Potential
under contract number HPMF-CT-2002-02031.

\section{Appendix}
\subsection{Dependence on a parameter} 
Let us consider a second order
differential equation of the following form
\begin{equation}
\label{eq:gsoe}
 y''= r(z,\varepsilon)y, \qquad '\equiv \frac{\rmd\phantom{z}}{\rmd z} ,
\end{equation} 
where $r(z,\varepsilon)$ is a rational function with respect to $z$ and
$\varepsilon$, i.e., $r\in\C(\varepsilon)(z)=\C(z,\varepsilon)$.  Here
$\varepsilon$ plays the role of a parameter. For a fixed value of
$\varepsilon$ we denote by $\cG^0(\varepsilon)$  the identity component of the
differential Galois group of equation~\eqref{eq:gsoe}. Let $U\subset \C$ denote an open not empty connected set with compact closure.    We show the
following.
\begin{lemma}
\label{lem:par}
Assume that:
\begin{enumerate}
\item Equation~\eqref{eq:gsoe} is Fuchsian.
\item For $\varepsilon\in U$,
   equation~\eqref{eq:gsoe} possesses
  $N$ singular points ($N$ does not depend on $\varepsilon$) for
  which exponents do not depend on $\varepsilon$.
\item For $\varepsilon_0\in U $,
  $\cG^0(\varepsilon_0)$ is not solvable (is not Abelian).
\end{enumerate}
Then, except finitely many values  of $\varepsilon\in U$, 
  $\cG^0(\varepsilon)$ is not solvable (is not Abelian).
\end{lemma}
\begin{proof}
  From the Kovacic algorithm in the form given
  in~\cite{Duval:92::,Morales:99::c} we know that, under our
  assumption, if $\cG^0(\varepsilon)$ is solvable, then there exists a
  polynomial $P$ (whose degree does not depend on $\varepsilon$) which
  is a solution of a linear differential equation $L(y)=0$ with
  coefficients in $\C(z,\varepsilon)$. The order of $L(y)=0$ does not
  depend on $\varepsilon$. We have a finite number of choices for the
  degree of $P$ and a finite number of choices of $L(y)=0$.  Finding a
  polynomial solution of linear equation $L(y)=0$ reduces to finding a
  non-trivial solution of a homogeneous linear system with
  coefficients in $\C(\varepsilon)$. But the last problem reduces to
  finding common zeros of a finite number of polynomials. We know that
  not all of these polynomials vanish identically (otherwise
  $\cG^0(\varepsilon_0)$ is solvable). Thus, there is at most a finite
  number of values of $ \varepsilon$ for which  they vanish
  simultaneously. Finally, let us notice that  set 
\[
\{\varepsilon\in U\,|\, \, \, 
  \cG^0(\varepsilon)\quad\text{is Abelian} \},
\]
is a subset of 
\[
\{\varepsilon\in U \,|\, \, \, 
  \cG^0(\varepsilon)\quad\text{is solvable} \}.
\]
\end{proof}
\subsection{Second order differential equations with rational coefficients}
Let us consider a second order
differential equation of the following form
\begin{equation}
\label{eq:gso}
 y''=r y, \qquad r\in\C(z), \qquad '\equiv \frac{\rmd\phantom{z}}{\rmd z} .
\end{equation} 
For this equation its differential Galois group $\cG$ is an algebraic
subgroup of $\mathrm{SL}(2,\C)$. The following lemma describes all
possible types of $\cG$ and relates these types to forms of solution
of \eqref{eq:gso}, see \cite{Kovacic:86::,Morales:99::c}.
\begin{lemma}
\label{lem:alg}
Let $\cG$ be the differential Galois group of equation~\eqref{eq:gso}.
Then one of four cases can occur.
\begin{enumerate}
\item $\cG$ is reducible (it is conjugated to a subgroup of triangular
  group) ; in this case equation \eqref{eq:gso} has an exponential
  solution of the form $y=\exp\int \omega$, where $\omega\in\C(z)$,
\item $\cG$ is conjugated with a subgroup of 
\[
D^\dag = \left\{ \begin{bmatrix} c & 0\\
                                0 & c^{-1}
                      \end{bmatrix}  \; \biggl| \; c\in\C^*\right\} \cup 
                      \left\{ \begin{bmatrix} 0 & c\\
                                c^{-1} & 0
                      \end{bmatrix}  \; \biggl| \; c\in\C^*\right\}, 
\]
  in this case equation
  \eqref{eq:gso} has a solution of the form $y=\exp\int \omega$, where
  $\omega$ is algebraic over $\C(z)$ of degree 2,
\item $\cG$ is primitive and finite; in this case all
  solutions of equation \eqref{eq:gso} are algebraic, 
  
\item $\cG= \mathrm{SL}(2,\C)$ and equation \eqref{eq:gso}
  has no Liouvillian solution.
\end{enumerate}
\end{lemma}

We  need  a more precise characterisation of case 1 in the above
lemma. It is given by the following lemma, see Lemma~4.2 in
\cite{Singer:93::a}.
\begin{lemma}
\label{lem:algc1}
  Let $\cG$ be the differential Galois group of
  equation~\eqref{eq:gso} and assume that $\cG$ is reducible.
  Then either
\begin{enumerate}
\item equation~\eqref{eq:gso} has a unique solution $y$ such that
  $y'/y\in\C(z)$, and $\cG$ is conjugate to a subgroup of the
  triangular group
\[
 \cT = \left\{ \begin{bmatrix} a & b\\
                        0& a^{-1}
               \end{bmatrix} \, |\, a,b\in\C, a\neq 0\right\}. 
\] 
Moreover, $\cG$ is a proper subgroup of $\cT$ if and only if there
exists $m\in\N$ such that $y^m\in\C(z)$. In this case $\cG$ is
conjugate to
\[
 \cT_m = \left\{ \begin{bmatrix} a & b\\
                        0& a^{-1}
               \end{bmatrix} \, |\, a,b\in\C, a^m=1\right\}, 
\]  
where $m$ is the smallest positive integer such that $y^m\in\C(z)$, or
\item equation~\eqref{eq:gso} has two linearly independent solutions
  $y_1$ and $y_2$ such that $y'_i/y_i\in\C(z)$, then $\cG$ is
  conjugate to a subgroup of
\[
 \cD = \left\{ \begin{bmatrix} a & 0\\
                        0& a^{-1}
               \end{bmatrix} \, |\, a\in\C, a\neq 0\right\}. 
\]
In this case, $y_1y_2\in\C(z)$. Furthermore, $\cG$ is conjugate to a
proper subgroup of $\cD$ if and only if $y_1^m\in\C(z)$ for some
$m\in\N$. In this case $\cG$ is a cyclic group of order $m$ where $m$
is the smallest positive integer such that $y_1^m\in\C(z)$.
\end{enumerate}
\end{lemma}

In case 2 of the above lemma we know that $v=y_1y_2\in\C(z)$.
Differentiating $v$ three times, and  using the fact that $y_i$ satisfies
equation \eqref{eq:gso},  we obtain
\begin{equation}
\label{eq:ssp}
 v'''= 2r' v + 4r v'.
\end{equation}
The above equation is called the second symmetric power of
equation~\eqref{eq:gso}. For applications of symmetric powers of
differential operators to study the existence of Liouvillian
solutions and differential Galois group see e.g.
\cite{Singer:93::a,Singer:95::,Ulmer:96::}.

To decide if case 2 from Lemma~\ref{lem:alg} occurs we can apply the
Kovacic algorithm. Here we present its part devoted to this case and
adopted to a Fuchsian equation.  At first we introduce notation.  We
write $r(z)\in\mathbb{C}(z)$ in the form
\begin{equation*}
\label{eq:rst}
r(z) = \frac{s(z)}{t(z)}, \qquad s(z),\, t(z) \in \mathbb{C}[z],
\end{equation*}
where $s(z)$ and $t(z)$ are relatively prime polynomials and $t(z)$ is
monic.  The roots of $t(z)$ are poles of $r(z)$. We denote $\Sigma':=
\{ c\in\mathbb{C}\,\vert\, t(c) =0 \}$ and
$\Sigma:=\Sigma'\cup\{\infty\}$.  The order $\ord(c)$ of $c\in\Sigma'$
is equal to the multiplicity of $c$ as a root of $t(z)$, the order of
infinity is defined by
\[ 
\mathrm{ord}(\infty):= \max(0, 4+\deg s - \deg t).
\]

Because  we assume that equation~\eqref{eq:gso} is Fuchsian,  we have
$\ord(c)\leq 2$ for $c\in\Sigma$. For each $c\in\Sigma'$ we have the
following expansion
\[
r(z) = \frac{a_c}{(z-c)^2} + O\left( \frac{1}{z-c}\right),
\]
and we define $\Delta_c = \sqrt{1+4a_c}$. For infinity we have

\begin{equation*}
r(z)=\dfrac{a_\infty}{z^2}+O\left(\dfrac{1}{z^3}\right),
\end{equation*}
and we define $\Delta_\infty = \sqrt{1+4a_\infty}$.

The algorithm consists of three steps.

\noindent
\textbf{Step I.}
For $c\in\Sigma'$ such that $\ord(c)=1$ we define 
 $E_c=\{4\}$;
if  $\ord(c)=2 $  
\[
 E_c := \left\{ 2 , 2(1+\Delta_c),
   2(1-\Delta_c)\right\}\cap\mathbb{Z}. 
\] 
If $\ord(\infty)<2$ we put $E_\infty=\{ 0,2,4\}$; if
$\ord(\infty)=2$ we define
\[
  E_\infty := \left\{2,  2(1+\Delta_\infty), 
  2(1-\Delta_\infty)\right\}\cap\mathbb{Z}. 
\] 
\textbf{Step II.} For each $e$ in the Cartesian product
\[
E:= {E}_{\infty}\times\prod_{c\in\Sigma'}{E}_c,
\]
we compute
\[
  d(e) := \frac{1}{2}\left(e_\infty- \sum_{c\in\Sigma'}e_c\right). 
\]
We select those elements $e\in E$ for which $d(e)$ is a
non-negative integer.  If there are no such elements Case 2 from
Lemma~\ref{lem:alg} cannot occur and the algorithm stops here.

\noindent 
\textbf{Step III.} For each element $e\in{E}$ such that
$d(e)=n\in\mathbb{N}_0$ we define
\[ 
\theta=\theta(z) = \frac{1}{2}
\sum_{c\in\Sigma'}\frac{e_c}{z-c},
\]
and we search for a monic polynomial $P=P(z)$ of degree $n$ satisfying the
following equation
\[ 
P''' + 3\theta P'' +(3 \theta^2 + 3\theta' -4r)P' + 
(\theta'' + 3 \theta\theta' + \theta^3 -4 r\theta - 2r')P =0.
\] 
If such polynomial exists, then equation~\eqref{eq:gso} possesses a
 solution of the form $w=\exp\int\omega$, where 
\[
\omega^2 + \psi\omega +\frac{1}{2}\psi' + \frac{1}{2}\psi^2 - r =0, \qquad 
\psi = \theta + \frac{P'}{P}.
\]
If we do not find such  polynomial, then case 2 in Lemma~\ref{lem:alg}
cannot occur.

\subsection{Riemann $P$ equation}

The Riemann $P$ equation \cite{Whittaker:35::} is the most general
second order differential equation with three regular singularities.
If we place, using homography, these singularities at $z=0,1,\infty$,
then it has the form
\begin{equation}
\label{eq:riemann}
\begin{split}
\dfrac{\mathrm{d}^2\xi}{\mathrm{d}z^2}&+\left(\dfrac{1-\alpha-\alpha'}{z}+
\dfrac{1-\gamma-\gamma'}{z-1}\right)\dfrac{\mathrm{d}\xi}{\mathrm{d}z}\\
&+
\left(\dfrac{\alpha\alpha'}{z^2}+\dfrac{\gamma\gamma'}{(z-1)^2}+
\dfrac{\beta\beta'-\alpha\alpha'-\gamma\gamma'}{z(z-1)}\right)\xi=0,
\end{split}
\end{equation}
where $(\alpha,\alpha')$, $(\gamma,\gamma')$ and $(\beta,\beta')$ are the
exponents at singular points. Exponents satisfy the Fuchs relation
\[
\alpha+\alpha'+\gamma+\gamma'+\beta+\beta'=1.
\]
We denote differences of exponents by
\[
\lambda=\alpha-\alpha',\qquad\nu=\gamma-\gamma',\qquad\mu=\beta-\beta'.
\]
For equation \eqref{eq:riemann} the necessary and sufficient
conditions for solvability of the identity component of its
differential Galois group are given by the following theorem due to
Kimura \cite{Kimura:69::}, see also \cite{Morales:99::c}.
\begin{theorem}[Kimura]
\label{thm:kimura}
  The identity component of the differential Galois group of
  equation~\eqref{eq:riemann} is solvable if and only if
\begin{itemize} 
\item[A:] at least one of  four numbers $\lambda+\mu+\nu$,
  $-\lambda+\mu+\nu$, $\lambda-\mu+\nu$, $\lambda+\mu-\nu$ is an odd
  integer, or
\item[B:] the numbers $\lambda$ or $-\lambda$ and $\mu$ or $-\mu$ and
  $\nu$ or $-\nu$ belong (in an arbitrary order) to some of the
  following fifteen families
\begin{center} 
\begin{tabular}{|c|c|c|c|c|} 
\hline 
1&$1/2+l$&$1/2+s$&arbitrary complex number&\\\hline 
2&$1/2+l$&$1/3+s$&$1/3+q$&\\\hline 
3&$2/3+l$&$1/3+s$&$1/3+q$&$l+s+q$ even\\\hline 
4&$1/2+l$&$1/3+s$&$1/4+q$&\\\hline 
5&$2/3+l$&$1/4+s$&$1/4+q$&$l+s+q$ even\\\hline 
6&$1/2+l$&$1/3+s$&$1/5+q$&\\\hline 
7&$2/5+l$&$1/3+s$&$1/3+q$&$l+s+q$ even\\\hline 
8&$2/3+l$&$1/5+s$&$1/5+q$&$l+s+q$ even\\\hline 
9&$1/2+l$&$2/5+s$&$1/5+q$&$l+s+q$ even\\\hline 
10&$3/5+l$&$1/3+s$&$1/5+q$&$l+s+q$ even\\\hline 
11&$2/5+l$&$2/5+s$&$2/5+q$&$l+s+q$ even\\\hline 
12&$2/3+l$&$1/3+s$&$1/5+q$&$l+s+q$ even\\\hline 
13&$4/5+l$&$1/5+s$&$1/5+q$&$l+s+q$ even\\\hline 
14&$1/2+l$&$2/5+s$&$1/3+q$&$l+s+q$ even\\\hline 
15&$3/5+l$&$2/5+s$&$1/3+q$&$l+s+q$ even\\\hline 
\end{tabular}\\[1.5ex] 
\end{center} 
Here $l,s$ and $q$ are integers. 
\end{itemize} 
\label{kimura} 
\end{theorem} 
The solvability conditions are sufficient for our purposes because if 
$\cG^0$ is not solvable, then obviously it is not Abelian.

\subsection{Lam\'e equation}

The Weierstrass 
form of the Lam\'e equation is following
\begin{equation} 
\label{eq:lame2} 
\frac{\mathrm{d}^2y}{\mathrm{d}t^2}=(\alpha\wp(t;g_2,g_3)+\beta)y, 
\end{equation} 
where $\alpha$ and $\beta$ are, in general, complex parameters and
$\wp(t;g_2,g_3)$ is the elliptic Weierstrass function with invariants
$g_2$, $g_3$. In other words, $\wp(t;g_2,g_3)$ is a solution of the
differential equation
\begin{equation*} 
\dot v^2=f(v),\qquad f(v)=4v^3-g_2v-g_3. 
\label{wei} 
\end{equation*} 
It is assumed that equation $f(v)=0$ has three different roots, so
\begin{equation*} 
\Delta=g_2^3 -27g_3^2 \neq  0.
\end{equation*} 
We recall that the modular function $j(g_2,g_3)$ associated with the
elliptic curve $u^2=4v^3-g_2v-g_3$ is defined as follows
\begin{equation*}
\label{eq:j}
j(g_2,g_3)=\frac{g_2^2}{g_2^3 -27g_3^2}.
\end{equation*}
Classically the Lam\'e equation is written with  parameter $n$
instead of $\alpha$ related by the formula $\alpha=n(n+1)$.  We see that
the Lam\'e equation depends on four parameters $(n,\beta,g_2,g_3)$.
The following lemma lists all the cases in which the identity
component of the differential Galois group of Lam\'e
equation~\eqref{eq:lame2} is Abelian, see
\cite[Sec.~2.8.4]{Morales:99::c}.
\begin{lemma} 
\label{lem:lame}
The identity component of the differential Galois group of Lam\'e
equation~\eqref{eq:lame2} is Abelian only in the following cases:
\begin{enumerate} 
\item the Lam\'e-Hermite case when $n\in\mathbb{Z}$ and three other
  parameters are arbitrary.
\item the Brioschi-Halphen-Crowford case for which
  $m:=n+\frac{1}{2}\in\mathbb{N}$, and remaining parameters
  $(g_2,g_3,\beta)$ satisfy an algebraic equation
\begin{equation*}
\label{eq:bdets}
Q_{m}\left(g_2,g_3,\beta\right)=0, 
\end{equation*}
see below.
\item the Baldassarri case $2n\not\in\Z$, and 
\[
n\pm q\in\Z\qquad\text{for some}\quad
  q\in\{1/4,1/6,1/10,3/10\},
\]
 with  additional algebraic restrictions
  on $(g_2,g_3,\beta)$. 
\end{enumerate} 
\end{lemma}

Polynomial $Q_m(g_2,g_3,\beta)$ which appears in the
Brioschi-Halphen-Crowford case, called the Brioschi determinant, is
defined as follows
\begin{equation}
\label{eq:bdet}
\left|\begin{array}{ccccccc}
\beta &m-1&0&0&0&\ldots& 0\\
q_{2,1}&\beta &2(m-2)&0&0&\ldots&0 \\
q_{3,1}&q_{3,2}&\beta&3(m-3)& 0& \dots&0\\
0 &q_{4,2}&q_{4,3}&\beta&4(m-4)&  \dots &0 \\
0 &0&\ddots&\ddots&\ddots&\ldots& 0\\
0 & 0& \ldots &                 q_{m-2,m-3}       & q_{m-1,m-2}    &\beta     &   m-1\\
0 & 0& \ldots & 0 &q_{m,m-2}&q_{m,m-1}&\beta
\end{array}
\right|,
\end{equation}
where
\[
q_{i+1,i} = \frac{g_2}{4}(2m-i)(m-i), \qquad 
q_{i+2,i} = \frac{g_3}{4}(2m-i)(2m-i-1).
\]

Algebraic restrictions on $(g_2,g_3,\beta)$ in the Baldassarri case
are involved. Instead of them we use the following proposition which
follows from one unpublished result of B.~Dwork,
see~\cite{Morales:99::c}.
\begin{proposition}
\label{prop:dwork}
The Baldassarri case for equation~\eqref{eq:lame2} occurs only for a finite number of values of $j(g_2,g_3)$.
\end{proposition}

%
%
\newcommand{\noopsort}[1]{}\def\cprime{$'$}
  \def\cydot{\leavevmode\raise.4ex\hbox{.}}

\end{document}